\documentclass[psamsfonts]{amsart}
\usepackage{graphicx}
\usepackage{amsmath,amssymb}
\usepackage{amsmath}
\usepackage{amsthm}
\usepackage{amssymb}
\usepackage{amscd}
\usepackage{amsfonts}
\usepackage{amsbsy}
\usepackage[noadjust]{cite} %The noadjust option makes (\cite{ref}) space correctly on the left.
\usepackage{epsfig,afterpage}
\usepackage{paralist}
\usepackage{psfrag}
\usepackage{mathrsfs}       %Gives the script font enabled by the invoked command \mathscr
\usepackage{verbatim}
\large\normalsize

\theoremstyle{plain}
\newtheorem{theorem}{Theorem}[section]

\theoremstyle{definition}

\newtheorem{example}[theorem]{Example}

\numberwithin{equation}{section}
\begin{document}
\newcommand{\T}{\mathbb{T}}
\newcommand{\R}{\mathbb{R}}
\newcommand{\Q}{\mathbb{Q}}
\newcommand{\N}{\mathbb{N}}
\newcommand{\Z}{\mathbb{Z}}
\newcommand{\tx}[1]{\quad\mbox{#1}\quad}
\newcommand{\yb}{{\bf y}}
\newcommand{\zb}{{\bf z}}
\newcommand{\vpb}{\pmb{\varphi}}

%Please put here any further newcommands.

%\noindent {\tt Neural, Parallel, and Scientific Computations 24 (2016) xx-xx}
\title[Volterra Equations]{\mbox{}\\[1cm]Symbolic Iterative Solution of Volterra Integral Equations}
\author{Hamid Semiyari}
\author{Douglas S. Shafer}
%\author{Author-3}
\begin{abstract}
\vspace{-.2cm}
\begin{center}
Mathematics Department, James Madison University, Harrisonburg, Virgina 22807,
USA
E-mail address: semiyahx@jmu.edu.

Mathematics Department, University of North Carolina at Charlotte, Charlotte,
North Carolina 28223, USA
E-mail address: dsshafer@uncc.edu
\end{center}
\vspace{-.2cm}
\end{abstract}
\thanks{\hspace{-.5cm}\tt
	Received xxxx, 2015
	\hfill 1061-5369 \$15.00 \copyright Dynamic Publishers, Inc.}
\maketitle
\thispagestyle{empty}

{\footnotesize
\noindent{\bf ABSTRACT.} In this work we show how auxiliary variables can be used to give an efficient and widely applicable method involving
symbolic manipulation and Picard iteration for approximating solutions of certain Volterra integral equations.

}

\section{Introduction and Preliminaries}\label{s1}
%Start the first section of the paper here.
Volterra integral equations of the second kind,
\begin{equation}\label{ve}
y(t) = \varphi(t) + \int_a^t K(t, s, y(s)) \, ds,
\end{equation}
lend themselves to solution by successive approximation using Picard iteration, although the process can break down
when quadratures that cannot be performed in closed form arise. In this article we offer a method for introducing
auxiliary variables in \eqref{ve} in the case that $K$ factors as $K(t, s, z) = f(t) k(s, z)$ in such a way that
\eqref{ve} embeds in a vector-valued polynomial Volterra integral equation, thus extending the method of auxiliary
variables, as expounded in \cite{PS} by Parker and Sochacki in the case of initial value problems, to the setting of
integral equations. We thereby obtain a computationally efficient method of \emph{symbolic} rather than numerical
computation for closely approximating solutions of \eqref{ve}. Of course the problem of impossible integrations could
also be addressed by replacing $\varphi$, $f$, and $k$ by initial segments of their power series expansions about $a$;
the method presented here seems to be an attractive alternative in some situations, among others those that involve 
denominators, like Example \ref{ex.ES} below, or those that involve powers of functions, like Example \ref{ex.BE}
below. It is equally easy to apply when the unknown function $y(t)$ appears in the argument of a transcendental
function, as in Example \ref{ex.trans}.

For reference we state the following generalization to the vector-valued case of Theorem 2.1.1 of \cite{H}. The 
proof in \cite{H} goes through with the obvious modifications. We only note that the proof is based on an application
of the Contraction Mapping Theorem. By way of notation, for a subset $S$ of $\R^m$ we let $C(S, \R^n)$ denoted the 
set of continuous mappings from $S$ into $\R^n$. 

\begin{theorem}\label{r:ve.exist.unique}
Let $I = [a, b] \subset \R$ and $J = \{ (x, y) : x \in I, \ y \in [a, x] \} \subset I \times I$. Suppose 
$\vpb \in C(I, \R^n)$ and $K \in C(J \times \R^n, \R^n)$ and that $K$ is Lipschitz in the last variable: there 
exists $L \in \R$ such that
\[
| K(x, y, \zb) - K(x, y, \zb') |_\textnormal{sum} \leqslant L | \zb - \zb' |_\textnormal{sum}
\]
for all $(x, y) \in J$ and all $\zb, \zb' \in \R^n$. Then the integral equation
\begin{equation}\label{e:vv.ve}
\yb(t) = \vpb(t) + \int_a^t K(t, s, \yb(s)) \, ds
\end{equation}
has a unique solution $\yb(t) \in C(I, \R^n)$.
\end{theorem}

Because the theorem was proved by means of the Contraction Mapping Theorem we immediately obtain the following result.

\begin{theorem}\label{r:ve.it.conv}
Under the hypotheses of Theorem \eqref{r:ve.exist.unique}, for any choice of the initial mapping $\yb^{[0]}(t)$ the 
sequence of Picard iterates
\[
\yb^{[k+1]}(t) = \vpb(t) + \int_a^t K(t, s, \yb^{[k]}(s)) \, ds
\]
converges to the unique solution of the integral equation \eqref{e:vv.ve}.
\end{theorem}

\section{The Method}\label{s2}
%Start the second section of the paper here.
Now let a Volterra integral equation
\begin{equation}\label{e:ve2}
y(t) = \varphi(t) + \int_a^t f(t) k(s, y(s)) \, ds
\end{equation}
be given, where $\varphi \in C([a,b], \R)$, $f \in C([a,b], \R)$, $k \in C([a,b] \times [a,b], \R)$, and $k$ satisfies
a Lipschitz condition in $y$. Introduce auxiliary variables $v_1, \dots , v_r$ in such a way that 
$\varphi = P(v_1, \dots, v_r)$, $f = Q(v_1, \dots, v_r)$, and $k  = R(y,v_1, \dots, v_r)$ (i.e.,  
$\varphi(t) = P(v_1(t), \dots, v_r(t))$, and so on), where $P$, $Q$, and $R$ are polynomials and the variables 
$v_1, \dots, v_r$ satisfy a system of first order polynomial ordinary differential equations
\begin{equation}\label{e:aux.ode}
\begin{aligned}
v_1' &= P_1(v_1, \dots, v_r) \\
       &\mspace{15mu}\vdots \\
v_r' &= P_r(v_1, \dots, v_r).
\end{aligned}
\end{equation}
To illustrate, suppose we wish to approximate the solution of $y(t) = 1 - \int_0^t \sin y(s) \, ds$. The integrand
indicates introducing $v_1 = v_1(t) = \sin y(t)$. Since $v_1'(t) = \cos t \, y'(t)$ we are then led to introduce
$v_2 = v_2(t) = \cos y(t)$, for which $v_2' = -v_1 \, y'$. From the integral equation itself we have
$y'(t) = -\sin y(t) = -v_1(t)$, so no additional auxiliary variables are needed; \eqref{e:aux.ode} is
$v_1' = -v_2 v_1$ and $v_2' = v_2^2$, and the integral equation is $y(t) = 1 - \int_0^t v_1(s) \, ds$. (In more
complicated situations some ingenuity can be required for this step. There are no known cases for which it has
proved impossible when the functions involved are analytic. See \cite{CPSW} for a fuller discussion.)

The initial value problem obtained by adjoining to \eqref{e:aux.ode} the initial conditions given by the values of 
$v_1$ through $v_r$ at $t = a$ has a unique solution, which is the unique solution of the vector-valued Volterra 
integral equation
\begin{equation}\label{e:aux.ve}
\begin{aligned}
v_1 &= v_1(a) + \int_a^t P_1(v_1(s), \dots, v_r(s)) \, ds \\
      &\mspace{15mu}\vdots \\
v_r &= v_r(a) + \int_a^t P_r(v_1(s), \dots, v_r(s)) \, ds,
\end{aligned}
\end{equation}
gotten simply by applying the Fundamental Theorem of Calculus to \eqref{e:aux.ode}. Adjoin to \eqref{e:aux.ve} the
original Volterra equation in the form
\[
y(t) = P(v_1(t), \dots, v_r(t))
      + \int_a^t Q(v_1(s), \dots, v_r(s)) R(y(s), v_1(s), \dots, v_r(s)) \, ds
\]
to obtain
\begin{equation}\label{e:aug.ve}
\begin{aligned}
y   &= P(v_1(t), \dots, v_r(t))
       + \int_a^t \, Q(v_1(s), \dots, v_r(s)) R(y(s), v_1(s), \dots, v_r(s)) \, ds \\
v_1 &= v_1(a) + \int_a^t P_1(v_1(s), \dots, v_r(s)) \, ds \\
    &\mspace{15mu}\vdots \\
v_r &= v_r(a) + \int_a^t P_1(v_1(s), \dots, v_r(s)) \, ds.
\end{aligned}
\end{equation}
System \eqref{e:aug.ve} satisfies the hypotheses of Theorem \eqref{r:ve.exist.unique}, hence has a unique solution, 
as does the original Volterra integral equation. Since $v_1, \dots, v_r$ are completely specified by 
\eqref{e:aux.ode} and \eqref{e:aux.ve}, the $y$ component of the solution of the augmented Volterra integral equation 
\eqref{e:aug.ve} must be the solution of \eqref{e:ve2}. But by Theorem \ref{r:ve.it.conv} the Picard iteration scheme 
applied to \eqref{e:aug.ve}, say with $y^{[0]}(t) \equiv \varphi(a)$ and $v_j^{[0]}(t) \equiv v_j(a)$, converges and 
is computationally feasible, so we obtain a computable approximation to the solution of \eqref{e:ve2}.

\section{Examples}\label{s3}
%Start the third section of the paper here.
In this section we illustrate the method by means of several examples which include both linear and nonlinear 
Volterra integral equations.

\begin{example}\label{ex.ES}
In \cite{Effati} Effati and Skandari introduced the linear Volterra integral equation of the second kind
\begin{equation}\label{Volt}
y(t) = e^t \sin t +\int_0^t \frac{2 + \cos t}{2 + \cos s} \, y(s) \, ds.
\end{equation}
The form of $\varphi(t)$ leads us to introduce $v_1 = e^t$ and $v_2 = \cos t$, and since $v_2' = - \sin t$,
also $v_3 = \sin t$. The integrand is then $(1 + v_2(t))(1+ v_2(s))^{-1} y(s)$; the denominator is the issue. To
express the integrand as a polynomial function of several variables, name the denominator $v_4 = 1 + v_2$ and its
reciprocal $v_5 = 1/ v_4$ so that the integral equation is 
$y(t) = v_1(t) v_2(t) + v_4(t) \int_0^t v_5(s) y(s) \, ds$. Taking the derivatives of the auxiliary variables introduced so far shows that no more are needed, so one appropriate choice of auxiliary variables is
\[
v_1 = e^t, \quad
v_2 = \cos t, \quad
v_3 = \sin t, \quad
v_4 = 2 + v_2, \quad
v_5 = \frac{1}{v_4},
\]
which satisfy the system of first order ordinary differential equations
\[
v_1' = v_1, \quad
v_2' = -v_3, \quad
v_3' = v_2, \quad
v_4' = v_2' = -v_3, \quad
v_5' = \frac{-v_4'}{v_4^2} = v_3v_5^2,
\]
which in turn is equivalent to
\begin{align*}
v_1(t) &= v_1(0) + \int_0^t v_1(s) \, ds\\
v_2(t) &= v_2(0) - \int_0^t v_3(s) \, ds\\
v_3(t) &= v_3(0) + \int_0^t v_2(s) \, ds\\
v_4(t) &= v_4(0) - \int_0^t v_3(s) \, ds\\
v_5(t) &= v_5(0) + \int_0^t v_3(s) v_5^2 \, ds.
\end{align*}
The initial values of the auxiliary variables are determined by their definition. The initial value $y(0)$ of the 
solution of the integral equation \eqref{Volt} is found simply by evaluating that equation at $t = 0$ to obtain 
$y(0) = 0$. Thus the iteration scheme is
\begin{align*}
y^{[k+1]}(t) &= v_1^{[k]}v_3^{[k]}+v_4^{[k]}\int_0^t\,v_5^{[k]}y^{[k]} ds\\
v_1^{[k+1]}(t) &= 1+\int_0^t\,v_1^{[k]} ds\\
v_2^{[k+1]}(t) &= 1-\int_0^t\,v_3^{[k]} ds\\
v_3^{[k+1]}(t) &= \int_0^t\,v_2^{[k]} ds\\
v_4^{[k+1]}(t) &= 3 - \int_0^t\,v_3^{[k]} ds\\
v_5^{[k+1]}(t) &= \frac{1}{3}+\int_0^t\,v_3^{[k]} (v_5^{[k]})^2 ds\\
\end{align*}
We can initialize as we please, but it is reasonable to choose
$y^{[0]}(t) \equiv y(0)$ and $v_j^{[0]}(t) \equiv v_j(0)$, i.e.,
$
(y^{[0]}, v_1^{[0]}, v_2^{[0]}, v_3^{[0]}, v_4^{[0]}, v_5^{[0]})(t) \equiv (0, 1, 1, 0, 3, \frac{1}{3})$. 

The exact solution of \eqref{Volt} is
\[
y(t) = e^t \sin t + e^t \Big(2 + \cos t \Big)\Big(\ln 3 - \ln\big(2 +\cos t\big)\Big),
\]
whose Maclaurin series, with its coefficients rounded to five decimal places, begins
\begin{align*}
y(t)       &= 1.00000 t + 1.50000 t^2 + 0.83333 t^3 + 0.16667 t^4 - 0.03333 t^5 \\
           &\mspace{282mu}- 0.02593 t^6 - 0.00529 t^7 + O(t^8). \\
\intertext{The Maclaurin series of the eighth Picard iterate, $y^{[8]}(t)$, with its coefficients rounded to 
five decimal places, begins}
y^{[8]}(t) &= 1.00000 t + 1.50000 t^2 + 0.83333 t^3 + 0.16667 t^4 - 0.03333 t^5 \\
           &\mspace{282mu}- 0.02593 t^6 - 0.00529 t^7 + O(t^8).
\end{align*}
The absolute value of the error in the approximation of the exact solution by $y^{[8]}(t)$ is practically zero up to
about $t = 0.4$, then increases monotonically to about 0.00057 at $t = 1$.
\end{example}

\begin{example}\label{ex.BE}
In \cite{Bia} Biazar and Eslami introduced the nonlinear Volterra integral equation of the second kind
\begin{equation}\label{NVolt8}
y(t) = \tfrac12 \sin 2t + \int_0^t \tfrac32 y(s)^2 \cos (s - t) \, ds.
\end{equation}
To fit this into the framework of \eqref{e:ve2} we begin by applying the cosine difference identity 
$\cos (t - s) = \cos s \cos t + \sin s \sin t$, obtaining
\[
y(t) = \tfrac12 \sin 2t
      + \tfrac32 \Big(\cos t \int_0^t y(s)^2 \cos s \, ds
                       + \sin t \int_0^t y(s)^2 \sin s \, ds\Big).
\]
Introducing the auxiliary variables $v = \cos t$ and $w = \sin t$, which solve the system
\[
v' = -w, \quad w' = v,
\]
upon integration we obtain the equivalent system of integral equations
\begin{align*}
v(t) &= v(0) - \int_0^t w(s) \, ds \\
w(t) &= w(0) + \int_0^t v(s) \, ds.
\end{align*}
The initial values of the auxiliary variables are determined by their definition. The initial value $y(0)$ of the 
solution of the integral equation \eqref{NVolt8} is found simply by evaluating that equation at $t = 0$ to obtain 
$y(0) = 0$. Thus the iteration scheme is
\begin{align*}
y^{[k+1]}(t)  &= w^{[k]} v^{[k]}
              + \tfrac32 \bigg(v^{[k]}(t) \int_0^t v^{[k]}(s) (y^{[k]})^2(s) \, ds
              + w^{[k]}(t) \int_0^t w^{[k]}(s) (y^{[k]})^2(s) \, ds\bigg) \\
w^{[k+1]}(t) &= 0 + \int_0^t v^{[k]}(s) \, ds \\
v^{[k+1]}(t) &= 1 - \int_0^t w^{[k]}(s) \, ds. \\
\end{align*}
We initialize with
\begin{align*}
y^{[0]}(t) &\equiv y(0) = 0 \\
w^{[0]}(t) &\equiv \sin 0 = 0 \\
v^{[0]}(t) &\equiv \cos 0 = 1.
\end{align*}
The exact solution of \eqref{NVolt8} is $y(t) = \sin t$, whose Maclaurin series, with its coefficients rounded to 
five decimal places, begins
\begin{align*}
y(t) = 1.00000 \, t - 0.16667 \, t^3 + 0.00833 \, t^5 - 0.00020 \, t^7 + O(t^9).
\end{align*}
The Maclaurin series of the eighth Picard iterate, $y^{[8]}(t)$, with its coefficients rounded to five decimal 
places, begins
\[
y^{[8]}(t) = 1.00000 \, t - 0.16667 \, t^3 + 0.008333 \, t^5 + 0.00000 \, t^7 + O(t^9).
\]
The absolute value of the error in the approximation of the exact solution by $y^{[8]}(t)$ is practically zero up to
about $t = 0.4$, then increases monotonically to about 0.001 at $t = 1$.
\end{example}

\begin{example}
As a somewhat more elaborate example consider the linear Volterra integral equation of the second kind 
given by
\begin{equation}\label{Volt9}
y(t) = \tan t - \tfrac14 \sin 2t - \tfrac12 t
       + \int_0^t \frac{1}{1+y^2(s)} \, ds.
\end{equation}
This is a corrected version of an integral equation given by Kamyad {\it et al} in \cite{Kamyad}. Because the 
integral part is independent of $t$, \eqref{Volt9} is equivalent to an initial value problem, namely
\[
y'(t) = \sec^2 t - \tfrac12 \cos 2t - \tfrac12
       + \frac1{1 + y^2(t)},
\quad
y(0) = 0.
\]
Of course by means of the identity $\cos^2 t = \frac12(1 + \cos 2t)$ the differential equation can be more 
compactly expressed as
\begin{equation}\label{e:yp.id}
y'(t) = \sec^2 t - \cos^2 t + \frac1{1 + y^2(t)},
\end{equation}
which will be important later.

To approximate the unique solution of \eqref{Volt9} we introduce the auxiliary variables
\[
v_1(t) =  \sin t, \quad
v_2(t) =  \cos t, \quad
v_3(t) =  \frac{1}{v_2}, \quad
v_4(t) =  1+y^2, \quad
v_5(t) = \frac{1}{v_4}.
\]

Note that in contrast with the previous examples the unknown function $y(t)$ figures into the definition of some of 
these variables, but in a polynomial way. Thus when we compute their derivatives $y$ also appears. Thanks to 
\eqref{e:yp.id}, it does so in a polynomial way, since by that identity $y' = v_3^2 - v_2^2 + v_5$ and we have 
additionally
\[
v_1' =  v_2, \quad
v_2' =  -v_1, \quad
v_3' =  v_1 v_3^2, \quad
v_4' =  2 y (v_3^2 - v_2^2 + v_5), \quad
v_5' =  -2 y v_5^2 (v_3^2 - v_2^2 + v_5).
\]

This system of ordinary differential equations, together with the equation satisfied by $y'$ and the known initial 
values of all the variables involved, is equivalent to the system of integral equations
\begin{align*}
y(t) &= v_1 v_3 - \tfrac12 v_1 v_2 - \tfrac12 t
            + \int_0^t v_5(s) \, ds \\
v_1(t) &=  \int_0^t v_2(s) \, ds \\
v_2(t) &=  1 -\int_0^t v_1(s) \, ds \\
v_3(t) &=  1 + \int_0^t v_1(s) v_3^2(s) \, ds  \\
v_4(t) &=  1 + 2 \int_0^t
                y(s) (v_3^2(s) - v_2^2(s) + v_5(s)) \, ds  \\
v_5(t) &=  1 - 2 \int_0^t
                y(s) v_5^2(s) (v_3^2(s) - v_2^2(s) + v_5(s)) \, ds.
\end{align*}
Setting up the obvious iteration scheme based on these integral equations, and initializing with the constant 
functions $y(t) \equiv y(0)$ and $v_j(t) \equiv v_j(0)$, the Picard iterate $y^{[28]}(t)$ with coefficients 
rounded to five decimal places is
\begin{align*}
y^{[28]}(t) &= 1.00000 \, t + 0.33333 \, t^3 + 0.13333 \,t^5 + 0.05397 \, t^7 + 0.02187 \, t^9 \\
            &\mspace{290mu}+ 0.00886 \, t^{11} + O(t^{13}).
\intertext{The exact solution is $y(t) = \tan t$, whose Maclaurin series, with coefficients rounded to five 
decimal places is}
y(t)        &= 1.00000 \, t + 0.33333 \, t^3 + 0.13333 \,t^5 + 0.05397 \, t^7 + 0.02187 \, t^9 \\
            &\mspace{290mu}+ 0.00886 \, t^{11} + O(t^{13}).
\end{align*}
On the interval $[0, 0.10]$ the error in the approximation of the exact solution by $y^{[28]}(t)$ increases
monotonically from zero to about $3.5 \times 10^{-14}$. 
\end{example}

Finally, we consider the Volterra equation already looked at in Section \ref{s2}, showing how the method
applies easily even when the unknown function is in the argument of a transcendental function.

\begin{example}\label{ex.trans}
The Volterra equation
\begin{equation}\label{ve.trans}
y = 1 - \int_0^t \, \sin y(s)\, ds
\end{equation}
has solution  $y(t) =2 \text{arccot}(\cot(\frac{1}{2})e^t)$ whose Maclaurin series with coefficients rounded to five
decimal places is
\begin{align*}
y(t) &= 1.00000 - 0.84147 \, t + 0.22732 \,t^2 + 0.05836 \, t^3 - 0.06154 \, t^4 + 0.00791 \, t^5 \\
     &\mspace{120mu}+ 0.01180 \, t^6 - 0.00629 \, t^7 - 0.00078 \,t^8 + 0.00202 \, t^9 + O \left(t^{10} \right).
\end{align*}
Introducing auxiliary variables as described for this example in Section \ref{s2} we obtain the recursion
\begin{align*}
y^{[k+1]}(t)   &= 1 - \int_0^t v_1^{[k]}(s) \, ds \\
v_1^{[k+1]}(t) &= \sin 1 - \int_0^t v_2^{[k]}(s) v_1^{[k]}(s) \, ds \\
v_2^{[k+1]}(t) &= \cos 1 + \int_0^t (v_1^{[k]})^2(s) \, ds.
\end{align*}
The ninth Picard iterate, $y^{[10]}(t)$, truncated to order nine and with its coefficients rounded to five decimal
places, is
\begin{align*}
y(t) &= 1.00000 - 0.84147 \, t + 0.22732 \,t^2 + 0.05836 \, t^3 - 0.06154 \, t^4 + 0.00791 \, t^5 \\
     &\mspace{120mu}+ 0.01180 \, t^6 - 0.00629 \, t^7 - 0.00078 \,t^8 + 0.00202 \, t^9 + O \left(t^{10} \right).
\end{align*}
The absolute value of the error in the approximation of the exact solution by $y^{[10]}(t)$ is practically zero up to
about $t = 0.5$, then increases monotonically to about $6.5 \times 10^{-4}$ at $t = 1$.
\end{example}

\section{Conclusion}\label{s4}
%Start the second section of the paper here.
After noting the extension to the vector-valued case of a well-known theorem on existence of solutions of Volterra
equations of the second kind, we have observed that the method of proof by means of the Contraction Mapping Theorem
guarantees that Picard iterates will converge to the solution. We have then described a method for introducing
auxiliary variables into Volterra equations of the form
\[
y(t) = \varphi(t) + \int_a^t f(t) k(s, y(s)) \, ds,
\]
in such a way that such an equation embeds in a vector-valued polynomial Volterra integral equation. We have thus
extended the method of auxiliary variables for surmounting the obstacle of impossible quadratures that can arise in
Picard iteration, well known in the case of initial value problems, to the setting of integral equations. We have
thereby obtained a computationally efficient method of symbolic rather than numerical computation for closely
approximating solutions of Volterra equations of this type, whether linear or nonlinear in the unknown solution $y$,
and even when $y$ appears in the argument of transcendental functions. We have illustrated the ease of use, broad
applicability, and efficiency of the method with examples.  

%\NPSCnewsection{ACKNOWLEDGMENTS}
%If there are any acknowledgments to research support etc., they go here. Notice that this section is unnumbered. 

\end{document}